\documentclass{colt2012} 

\usepackage{algorithm,algorithmic}
\usepackage{amssymb, multirow, paralist, color}
\newtheorem{thm}{Theorem}

\newtheorem{ass}[thm]{Assumption}
\usepackage{enumitem}
\usepackage{textcase}
\usepackage[T1]{fontenc}

\def \R {\mathbb{R}}

\def \x {\mathbf{x}}

\def \E {\mathrm{E}}

\def \x {\mathbf{x}}

\def \z {\mathbf{z}}

\def \y {\mathbf{y}}
\def \u {\mathbf{u}}

\def \g {\mathbf{g}}

\def \y {\mathbf{y}}
\def \E {\mathrm{E}}
\def \x {\mathbf{x}}
\def \g {\mathbf{g}}

\def \z {\mathbf{z}}
\def \u {\mathbf{u}}

\def \R {\mathbb{R}}

\def \N {\mathcal{N}}

\def \G {\mathcal {G}}

\begin{document}

\title[Non-convex Extragradient]{On the Convergence of (Stochastic) Gradient Descent with Extrapolation for Non-Convex Optimization}
\author{\Name{Yi Xu}$^\dagger$\Email{yi-xu@uiowa.edu}\\
\Name{Zhuoning Yuan}$^\dagger$\Email{zhuoning-yuan@uiowa.edu}\\
 \Name{Sen Yang}$^\ddagger$ \Email{senyang.sy@alibaba-inc.com}\\
 \Name{Rong Jin}$^\ddagger$\Email{jinrong.jr@alibaba-inc.com}\\
\Name{Tianbao Yang}$^\dagger$\Email{tianbao-yang@uiowa.edu}\\
   \addr $^\dagger$Department of Computer Science, The University of Iowa, Iowa City, IA 52242, USA  \\
   \addr$^\ddagger$Machine Intelligence Technology, Alibaba Group, Bellevue, WA 98004, USA
}
\maketitle
\vspace*{-0.5in}
\begin{center}
First version: January 29, 2019
\end{center}

\begin{abstract}
Extrapolation is a well-known technique for solving convex optimization and variational inequalities and recently attracts some attention for non-convex optimization. Several recent works have empirically shown its success in some machine learning tasks. However, it has not been analyzed for non-convex minimization and there still remains a gap between the theory and the practice. In this paper, we analyze gradient descent and stochastic gradient descent with extrapolation for finding an approximate first-order stationary point in smooth non-convex optimization problems. Our convergence upper bounds show that the algorithms with extrapolation can be accelerated than without extrapolation.  
\end{abstract}
\section{Introduction}
We are interested in solving the following {\bf non-convex} optimization problem:
\begin{align}\label{opt:extra:grad:1}
\min_{\x\in\R^d}f(\x),
\end{align}
where $f(\x)$ is $L$-smooth. When the objective function can be written as an expectation of a random function, then (\ref{opt:extra:grad:1}) becomes a {\bf stochastic non-convex} optimization problem:
\begin{align}\label{opt:extra:grad:2}
\min_{\x\in\R^d}f(\x) := \E[f(\x;\xi)],
\end{align}
where $\xi$ is a random variable. In this paper, we consider solving the problem~(\ref{opt:extra:grad:1}) by gradient descent with extrapolation (GDE) method and solving the problem~(\ref{opt:extra:grad:2}) by stochastic gradient descent with extrapolation (SGDE) method. Non-convex optimization has brought tremendous success in many areas of machine learning including deep learning~\citep{goodfellow2016deep}, tensor decomposition~\citep{ge2015escaping}, and low-rank matrix completion~\citep{jain2013low}. Many existing works have shown non-convex losses may yield improved robustness and classification accuracy~\citep{chapelle2009tighter, wu2007robust, nguyen2013algorithms,DBLP:journals/corr/abs-1805-07880}. It is well known that traditional gradient descent (GD) and its variants such as stochastic gradient descent (SGD) are widely used in solving the problems~(\ref{opt:extra:grad:1}) and (\ref{opt:extra:grad:2}), respectively. The convergence results are also well-studied both for GD and SGD methods~\citep{nesterov1998introductory, ghadimi2013stochastic,yangnonconvexmo}. For example, \cite{nesterov1998introductory} has shown that GD enjoys iteration complexity of $O(1/\epsilon^2)$ for finding an $\epsilon$-first-order stationary point (i.e., find an $\x$ such that $\|\nabla f(\x)\| \leq \epsilon$) of problem~(\ref{opt:extra:grad:1}). \cite{ghadimi2013stochastic} established a $O(1/\epsilon^4)$ iteration complexity of SGD for finding an $\epsilon$-first-order stationary point in expectation satisfying that $\E[\|\nabla f(\x)\|] \leq \epsilon$ for solving the problem~(\ref{opt:extra:grad:2}). \cite{yangnonconvexmo} then extended the result to stochastic momentum methods.
Although GD and SGD achieve lots of success, recent works have shown that extragradient descent methods perform better or converge faster than GD/SGD in several machine learning tasks such as training generative adversarial networks (GANs)~\citep{yadav2017stabilizing,gidel2018variational}, training low bit neural network~\citep{leng2018extremely}, learning Gaussian mixture models~\citep{mertikopoulos2018mirror}, and solving linear programming~\citep{wang2009large}. However, the theoretical guarantee of non-asymptotical convergence of GDE (resp. SGDE) is still unclear for the general non-convex problem (\ref{opt:extra:grad:1}) (resp. problem (\ref{opt:extra:grad:2})). 
In this paper, we analyze  GDE and a new variant of SGDE and establish their convergence results for finding an approximate first-order stationary point in non-convex optimization.
The main contributions of this paper are summarized as follows.
\begin{itemize}
\item We analyze a variant of GDE for a general smooth non-convex problem~(\ref{opt:extra:grad:1}).  It enjoys an iteration complexity of $O(1/\epsilon^2)$ for finding an $\epsilon$-first-order stationary point $\x$ of problem (\ref{opt:extra:grad:1}) that satisfying $\|\nabla f(\x)\| \leq \epsilon$. Our convergence bound shows that it could be faster than the GD method.  
\item We analyze mini-batch SGDE algorithm, showing that it achieves a total gradient complexity of $O(1/\epsilon^4)$ for finding an $\epsilon$-first-order stationary solution $\x$ of problem (\ref{opt:extra:grad:2}) in expectation with a mini-batch size of $O(1/\epsilon^2)$. To avoid the mini-batch requirement, we also propose a stagewise SGDE, which enjoys the same gradient complexity of $O(1/\epsilon^4)$ but without the requirement of a large mini-batch size. Our convergence bound also shows that it could achieve practical speed-up. 

\end{itemize}

\section{Related Work}
In this section, we review some related work about gradient descent with extrapolation methods. 
The extragradient method was first introduced by~\cite{korpelevich1976extragradient, korpelevichextrapolation} for solving variational inequality problems (VIP)~\citep{hartman1966some,harker1990finite}, i.e., finding a point $\x_*\in\Omega$ such that
$\langle \G(\x), \x_* - \x\rangle\le 0, \forall \x\in\Omega$,
where $\Omega$ is a nonempty closed convex subset of $\R^d$ and $\G:\R^d \to\R^d$ is an operator.
It generates a pair of sequence by carrying out two projections in each iteration:
\begin{align*}
\x_t &= \text{P}_\Omega[\z_{t-1} - \eta \G(\z_{t-1})],\\
\z_t &= \text{P}_\Omega[\z_{t-1} - \eta \G(\x_t)].
\end{align*}
Most subsequent research works~e.g., \citealt{tseng2000modified, censor2011subgradient, thong2018weak,nemirovski2004prox,nesterov2007dual,solodov1999hybrid,monteiro2010complexity, monteiro2011complexity,juditsky2011solving,chen2017accelerated} 
have analyzed the convergence of extragradient method and its variants for solving (stochastic) VIP under the assumptions of $L$-Lipschtiz continuous and monotone operator $\G$.  If one considers the minimization of a function as a VI problem, then Lipschitz continuous and monotone operator means  the gradient of a convex function that is Lipschitz continuous. It is also notable that~\cite{gidel2018variational} proposed stochastic extragradient algorithms for solving min-max saddle point problems from a perspective of variational inequality. In the theoretical side, their convergence rates of proposed algorithms are built based on an assumption that the considered problems are convex-concave or the variational inequalities are monotone.
Few works have considered (stochastic) extragradient methods for non-monotone VI~\citep{kannan2014optimal,DBLP:journals/coap/DangL15} under some pseudo-monotonicity assumption. 
In contrast, we directly analyze GDE methods and their convergence for finding a stationary point of  smooth non-convex optimization problems without considering the above assumptions.

In the context of optimization setting, extragradient method and its accelerated/extended version were well studied with  the establishments of convergence rate.  
It has been shown~\citep{luo1993error, wang2014iteration} that extragradient method is a special case of feasible descent method (FDM).  Under local error bound assumption, \cite{luo1993error} have proved linear convergence of extragradient method for solving convex optimization problems. \cite{monteiro2013accelerated} applied  hybrid proximal extragradient (HPE) method to convex optimization by proposing an accelerated HPE, enjoying the convergence rate of $O(1/T^2)$.
Recently, \cite{diakonikolas2018accelerated} developed an accelerated extragradient descent (AXGD) method for solving smooth and convex problems by combining the key ideas from Nesterov's accelerated gradient (NAG) method~\citep{nesterov1983method} and Nemirovski's mirror-prox method~\citep{nemirovski2004prox}. AXGD achieved a convergence rate of $O(1/T^2)$, matching the order of NAG's convergence rate. 
\cite{DBLP:journals/jmlr/ChiangYLMLJZ12,DBLP:journals/ml/YangMJZ14} and several subsequent works~\citep{DBLP:conf/nips/RakhlinS13,DBLP:conf/colt/RakhlinS13} have considered the extragradient method for online convex optimization that repeatedly use an online gradient for two updates,  and showed smaller regret compared with online gradient method for smooth functions.  

Very recently, \cite{nguyen2018extragradient} proposed an extended extragradient method (EEG) to minimize the sum of two functions that one is smooth and another is convex. EEG uses two proximal gradient steps at each iteration, which is slightly different from two orthogonal projection steps of classical extragradient. Like classical extragradient method, EEG still has the issue that computing two gradients might seriously affect the efficiency of the algorithm. For non-convex case, under the Kurdyka-\L ojasiewicz (KL) assumption~\citep{bolte2007lojasiewicz, bolte2010characterizations, bolte2017error}, they have shown that the sequence generating by EEG converges to a first-order critical point of the considered problem with finite length. Their convergence rate is asymptotic and heavily depends on the \L ojasiewicz exponent parameter $\theta$~\citep{bolte2017error}, which value is specific to the particular problem. 
By contrast, we consider GDE methods for solving general smooth but non-convex problems, and estbalish  a non-asymptotic convergence result with an iteration complexity of $O(1/\epsilon^2)$ for finding an $\epsilon$-first-order stationary point with potential improvement than the GD method. We also propose two variants of GDE method in stochastic setting, namely mini-batch SGDE and stagewise SGDE with both of them achieving an iteration complexity of $O(1/\epsilon^4)$ for finding an $\epsilon$-first-order stationary point in expectation. 
It is worth mentioning that our GDE and SGDE methods only need to compute gradient or stochastic gradient once per iteration inspired by~\citep{DBLP:journals/jmlr/ChiangYLMLJZ12,DBLP:journals/ml/YangMJZ14}, implies that our methods are more efficient than EEG since it saves the computation of (stochastic) gradient in each updating.


\section{Preliminaries}
In this section, we will present some notations and a previous result about extrapolation. Recall the problem of interest is 
\begin{align}\label{opt:extra:grad:3}
\min_{\x\in\R^d}f(\x),
\end{align}
or
\begin{align}\label{opt:extra:grad:4}
\min_{\x\in\R^d}f(\x) := \E[f(\x;\xi)],
\end{align}
where $\xi$ is a random variable, both $f(\x)$ and $f(\x;\xi)$ are non-convex functions. Let us denote by $\x_*$ the global minimum of $f(\x)$, i.e., $\x_* \in \arg\min_{\x\in\R^d} f(\x)$.
First, we make the following assumptions throughout the paper.
\begin{ass}\label{ass:1} 
\begin{itemize} [leftmargin=0.35in]
\item[(i).] $f(\x)$ has $L$-Lipschitz continuous gradient, i.e., there exists $L>0$ such that $\| \nabla f(\x) - \nabla f(\y)\| \leq L \|\x-\y\| $ for all $\x,\y\in\R^d$;
\item[(ii).] there exists $\Delta<\infty$ such that $f(\x) - f(\x_*)\leq \Delta$ for all $\x\in\R^d$;
\item[(iii).] every random function $f(\x; \xi)$ is differentiable;
\item[(iv).] there exists $G>0$ such that $\mathbb{E}[\|\nabla f(\x; \xi) - \nabla f(\x)\|^2 ]\leq G^2$ holds.
\end{itemize}
\end{ass}
{\bf Remark. }Assumption~\ref{ass:1} (iii) (iv) are standard assumptions made in the literature of stochastic non-convex optimization~\citep{ghadimi2013stochastic,yangnonconvexmo}. 
Assumption~\ref{ass:1} (ii) is used to get the iteration complexity of an algorithm, in particular, it is needed in getting the iteration complexity of Stagewise SGDE. While for iteration complexities of GDE and SGDE, we only need a weaker assumption that for an initial solution $\x_0\in\R^d$, there exists a constant $\Delta_0 >0$ such that $f(\x_0)-f(\x_*)\leq \Delta_0$.

Next, to measure the convergence of non-convex and smooth optimization problems as in\citep{nesterov1998introductory, ghadimi2013stochastic,yangnonconvexmo}, we need the following definition of first-order stationary point.
\begin{definition}[First-order stationary point] For problem (\ref{opt:extra:grad:3}) or (\ref{opt:extra:grad:4}), a point $\x\in\R^d$ is called a first-order stationary
point if
\begin{align*}
\|\nabla f(\x)\| = 0.
\end{align*}
Moreover, if
\begin{align*}
\|\nabla f(\x)\| \leq \epsilon,
\end{align*}
then the point $\x$ is said to be an $\epsilon$-stationary point.
\end{definition}

We then introduce the Moreau envelope function of $f(\x)$ and proximal mapping, which are formally stated as follows.
\begin{definition}
For any $\lambda > 0$, the following function is called a Moreau envelope of $f$
\begin{align}
f_{\gamma}(\x) := \min_{\y\in\R^d}\left\{f(\y) + \frac{1}{2\gamma}\|\y-\x\|^2 \right\}.
\end{align}
Moreover, the optimal solution to the above problem is called a proximal mapping of $f$:
\begin{align}
\text{prox}_{\gamma f}(\x) := \arg \min_{\y\in\R^d}\left\{f(\y) + \frac{1}{2\gamma}\|\y-\x\|^2 \right\}\end{align}
\end{definition}
Let $\widehat\x := \text{prox}_{\gamma f}(\x)$, 
it has been shown that~\citep{davis2018stochastic}
\begin{align}
\nabla f_{\gamma}(\x) = \frac{1}{\gamma}(\x-\widehat\x)
\end{align}
and 
\begin{align}\label{moreau:ineqs}
f(\widehat \x) \leq f(\x),~~\|\x-\widehat\x\| = \gamma\|\nabla f_{\gamma}(\x) \|,~~\| \nabla f(\widehat \x)\| \leq \| \nabla f_{\gamma}(\x) \|.
\end{align}

Finally, we will end up with a key lemma in~\citep{nemirovski2004prox} for our analysis.
\begin{lemma}[\textbf{Lemma 3.1},~\cite{nemirovski2004prox}]\label{lem:key1}
Let $\omega(\z)$ be a  $\alpha$-strongly convex function with respect to the norm $\|\cdot\|$, whose dual norm is denoted by $\|\cdot\|_*$,  and $D(\x,\z) = \omega(\x)- (\omega(\z) + (\x-\z)^{\top}\omega'(\z))$ be the Bregman distance induced by function $\omega(\x)$. Let $Z$ be a convex compact set, and $U\subseteq Z$ be convex and closed.  Let $\z\in Z$, $\gamma>0$, Consider the points,
\begin{align}
\x &= \arg\min_{\u\in U} \gamma\u^{\top}\xi + D(\u, \z)\label{eqn:project1},\\
\z_+&=\arg\min_{\u\in U} \gamma\u^{\top}\zeta + D(\u,\z),\label{eqn:project2}
\end{align}
then for any $\u\in U$, we have
\begin{align}\label{eqn:ineq}
\gamma\zeta^{\top}(\x-\u)\leq  D(\u,\z) - D(\u, \z_+) + \frac{\gamma^2}{\alpha}\|\xi-\zeta\|_*^2 - \frac{\alpha}{2}[\|\x-\z\|^2 + \|\x-\z_+\|^2].
\end{align}
\end{lemma}

\section{Main Results}
In this section, we will present the proposed algorithms and the main results of their convergence. We will first introduce a GDE algorithm for solving the problem~(\ref{opt:extra:grad:3})  and a mini-batch SGDE algorithm for solving the problem~(\ref{opt:extra:grad:4}). Then we will extend the mini-batch SGDE algorithm to stagewise SGDE without using a mini-batch of samples, which is more practical and user-friendly. 

\subsection{Gradient Descent with Extrapolation} 
The detailed updating steps of GDE are described in Algorithm~\ref{alg:GDE}, where $\eta>0$ is the step size. 
Please note that the updates of our GDE is slightly different from the updates of traditional GDE:
$\x_t = \z_{t-1} - \eta \nabla f(\z_{t-1}),~\z_t = \z_{t-1} - \eta \nabla f(\x_t).$
One issue of the traditional GDE is that the algorithm alternately computes the gradients at two points $\{\z_t\}$ and $\{\x_t\}$ for each iteration, implying that it is twice costly than the GD method that computes gradient per-iteration. By contrast, our considered GDE method stores and reuses the previous gradient to update the new extrapolation point. That is to say, our GDE only requires computing gradient once per-iteration. The similar idea was used in the online convex optimization~\citep{DBLP:journals/ml/YangMJZ14,DBLP:journals/jmlr/ChiangYLMLJZ12} and recently by~\cite{yadav2017stabilizing,gidel2018variational} for training GAN. 
In this paper, we focus on analyzing the convergence of GDE for non-convex optimization, and the result is presented in Theorem~\ref{thm:GDE}.
\begin{algorithm}[t]
\caption{GDE} \label{alg:GDE}
\begin{algorithmic}
\STATE Initilization: $\z_0 = \x_0$, $\g_0 = \nabla f(\x_{0})$
\FOR{$t=1,\ldots, T$}
    	\STATE $\x_t = \z_{t-1} - \eta \g_{t-1} $
	\STATE $\g_t = \nabla f(\x_{t})$
	\STATE $\z_t = \z_{t-1} - \eta \g_t$
\ENDFOR
\end{algorithmic}
\end{algorithm}
\begin{algorithm}[t]
\caption{Mini-batch SGDE} \label{alg:SGDEmini}
\begin{algorithmic}
\STATE Initilization: $\z_0 = \x_0$ and $\g_0 = \frac{1}{m}\sum_{i=1}^{m}\nabla f(\x_{0};\xi_{i,0})$
\FOR{$t=1,\ldots, T$}
    	\STATE $\x_t = \z_{t-1} -   \eta \g_{t-1}$
	\STATE $\g_t =  \frac{1}{m}\sum_{i=1}^{m}\nabla f(\x_t;\xi_{i,t}) $
	\STATE $\z_t =  \z_{t-1} -  \eta \g_t$
\ENDFOR
\end{algorithmic}
\end{algorithm}

\begin{thm}\label{thm:GDE}
Under Assumption~\ref{ass:1} (i), let $\eta \leq \frac{1}{12L}$ and $\x_1=\z_0= \x_0$, then GDE ensures that
\begin{align}\label{upp:bound:GDE}
\min_{t\in\{1,\ldots, T\}}\|\nabla f(\x_t)\|^2\leq 
 \frac{8(f(\x_0) - f(\x_*))}{\eta T} - \frac{1}{\eta^2T}\sum_{t=0}^{T-1}\|\x_{t+1} -\x_{t}\|^2,
\end{align}
where $\x_* = \arg\min_{\x\in\R^d}f(\x)$. Under an additional assumption that $f(\x_0) - f(\x_*) \leq \Delta_0$ where $\Delta_0>0$ is a constant, in particular in order to have $\min_{t\in\{1,\ldots, T\}}\|\nabla f(\x_t)\| \leq \epsilon$, the iteration complexity is $T=O(1/\epsilon^2)$.
\end{thm}
{\bf Remark. }The iteration complexity $O(1/\epsilon^2)$ of GDE is at least the same order of the GD method for smooth non-convex optimization. However, comparing with the convergence upper bound of GD, the above bound of GDE in (\ref{upp:bound:GDE}) has an additional negative term $- \frac{1}{\eta^2T}\sum_{t=0}^{T-1}\|\x_{t+1} -\x_{t}\|^2$, which should be beneficial for accelearting convergence in practice. 
\begin{proof}
By the $L$-smooth of $f(\x)$ we have
\begin{align*}
f(\x_t)&\leq f(\z_{t-1}) + \nabla f(\z_{t-1})^{\top}(\x_t - \z_{t-1}) + \frac{L}{2}\|\x_t - \z_{t-1}\|^2\\
&\leq f(\z_{t-1}) + \nabla f(\x_t)^{\top}(\x_t - \z_{t-1}) + (\nabla f(\z_{t-1}) - \nabla f(\x_t))^{\top}(\x_t - \z_{t-1}) + \frac{L}{2}\|\x_t - \z_{t-1}\|^2\\
&\leq f(\z_{t-1}) + \nabla f(\x_t)^{\top}(\x_t - \z_{t-1}) + \frac{3 L}{2}\|\x_t - \z_{t-1}\|^2.
\end{align*}
Applying Lemma~\ref{lem:key1} with $\u=\z_{t-1}, \x=\x_t, \z=\z_{t-1}, \z_+=\z_t, \xi=\nabla f(\x_{t-1}), \zeta=\nabla f(\x_t), \gamma=\eta$, we have
\begin{align*}
&\nabla f(\x_t)^{\top}(\x_t - \z_{t-1})\\
\leq& -\frac{\|\z_t - \z_{t-1}\|^2}{2\eta} + \eta \|\nabla f(\x_t) - \nabla f(\x_{t-1})\|^2 - \frac{1}{2\eta}(\|\x_t - \z_{t-1}\|^2 + \|\x_t - \z_t\|^2).
\end{align*}
Combining the above two inequalities together, we have
\begin{align*}
f(\x_t)\leq& f(\z_{t-1}) -\frac{\|\z_t - \z_{t-1}\|^2}{2\eta} + \eta \|\nabla f(\x_t) - \nabla f(\x_{t-1})\|^2 \\
&- \frac{1}{2\eta}(\|\x_t - \z_{t-1}\|^2 + \|\x_t - \z_t\|^2) + \frac{3L}{2}\|\x_t -\z_{t-1}\|^2.
\end{align*}
Moreover, by the smoothness of $f(\x)$
\begin{align*}
f(\z_t)&\leq f(\x_t) + \nabla f(\x_t)^{\top}(\z_t - \x_t) + \frac{L}{2}\|\x_t -\z_t\|^2\\
&\leq f(\x_t) + \frac{\eta\|\nabla f(\x_t)\|^2 }{4} + \frac{\|\x_t - \z_t\|^2}{\eta} + \frac{L}{2}\|\x_t - \z_t\|^2.
\end{align*}
Combining the above two inequalities together, we have
\begin{align*}
f(\z_t)\leq& f(\z_{t-1})+ \left(\frac{\eta\|\nabla f(\x_t)\|^2 }{4}   -\frac{\|\z_t - \z_{t-1}\|^2}{2\eta}\right) +  \bigg(\frac{1}{\eta}+\frac{L}{2}\bigg)\|\x_t - \z_t\|^2  \\
&+ \eta \|\nabla f(\x_t) - \nabla f(\x_{t-1})\|^2- \frac{1}{2\eta}(\|\x_t - \z_{t-1}\|^2 + \|\x_t - \z_t\|^2) + \frac{3L}{2}\|\x_t -\z_{t-1}\|^2\\
= & f(\z_{t-1})- \frac{\eta\|\nabla f(\x_t)\|^2 }{4}  +  \bigg(\eta+\frac{L\eta^2}{2}+\eta\bigg)\|\nabla f(\x_t) - \nabla f(\x_{t-1})\|^2 \\
&- \frac{1}{2\eta}(\|\x_t - \z_{t-1}\|^2 + \|\x_t - \z_t\|^2) + \frac{3\eta^2 L}{2}\|\nabla f(\x_{t-1})\|^2,
\end{align*}
where we use the facts that $\z_t - \z_{t-1} = - \eta \nabla f(\x_t)$, $\z_t - \x_t =  - \eta(\nabla f(\x_t) - \nabla f(\x_{t-1}))$ and $\x_t - \z_{t-1} = - \eta \nabla f(\x_{t-1})$. 
Taking summation on both sides, we have 
\begin{align*}
&\sum_{t=1}^T\frac{\eta}{4}\|\nabla f(\x_t)\|^2 - \frac{3L\eta^2}{2}\sum_{t=1}^T\|\nabla f(\x_{t-1})\|^2\\
\leq& (\eta + \frac{\eta^2 L}{2}  + \eta)\sum_{t=0}^{T-1}\|\nabla f(\x_{t+1}) - \nabla f(\x_{t})\|^2\\
& - \frac{1}{2\eta}\sum_{t=1}^T(\|\x_t - \z_{t-1}\|^2  + \|\x_t - \z_t\|^2) +\sum_{t=1}^T \bigg(f(\z_{t-1}) - f(\z_t)\bigg)
\end{align*}
Note that 
\begin{align*}
&\sum_{t=1}^T(\|\x_t - \z_{t-1}\|^2  + \|\x_t - \z_t\|^2)\\
=&\sum_{t=0}^{T-1}\|\x_{t+1} - \z_{t}\|^2  + \sum_{t=1}^T\|\x_t - \z_t\|^2\\
=&\sum_{t=1}^{T-1}\|\x_{t+1} - \z_{t}\|^2  + \|\x_{1} - \z_{0}\|^2 + \sum_{t=1}^T\|\x_t - \z_t\|^2 \\
\geq& \frac{1}{2}\sum_{t=1}^{T-1}\|\x_t - \x_{t+1}\|^2+ \|\x_{1} - \z_{0}\|^2\\
=  &\frac{1}{2}\sum_{t=0}^{T-1}\|\x_t - \x_{t+1}\|^2 + \frac{1}{2}\|\x_1 - \x_{0}\|^2.~~(\text{Since }\z_0 = \x_0)
\end{align*}
Combining the above inequalities together, we have
\begin{align*}
\sum_{t=1}^T\bigg(\frac{\eta}{4} - \frac{3L\eta^2}{2}\bigg)\|\nabla f(\x_t)\|^2\leq& \frac{3L\eta^2}{2}\|\nabla f(\x_0)\|^2 -  \frac{1}{4\eta}\|\x_1 - \x_{0}\|^2+\sum_{t=1}^T \bigg(f(\z_{t-1}) - f(\z_t)\bigg)\\
&+(2\eta L^2 + \frac{\eta^2 L^3}{2} - \frac{1}{4\eta})\sum_{t=0}^{T-1}\|\x_{t+1} -\x_{t}\|^2.
\end{align*}
Since $\eta \leq \frac{1}{12L}$, then $2\eta L^2  + \frac{\eta^2 L^3}{2} - \frac{1}{4\eta} \leq  -\frac{1}{8\eta}$ and $\frac{\eta}{4} - \frac{3L\eta^2}{2} \geq \frac{\eta}{8}$.
Note that we can define $\nabla f(\x_0)=0$ (i.e, $\x_1 = \z_0$), which will not affect our analysis above. Then $\x_1=\z_0= \x_0$,  and 
\begin{align*}
\frac{\eta}{8}\sum_{t=1}^T\|\nabla f(\x_t)\|^2&\leq f(\z_0) - f(\x_*) - \frac{1}{8\eta}\sum_{t=0}^{T-1}\|\x_{t+1} -\x_{t}\|^2 ,
\end{align*}
which implies
\begin{align*}
\min_{t\in\{1,\ldots, T\}}\|\nabla f(\x_t)\|^2 \le \frac{1}{T}\sum_{t=1}^T\|\nabla f(\x_t)\|^2 \leq \frac{8(f(\z_0) - f(\x_*))}{\eta T} - \frac{1}{\eta^2T}\sum_{t=0}^{T-1}\|\x_{t+1} -\x_{t}\|^2.
\end{align*}
\end{proof}
\subsection{Stochastic Gradient Descent with Extrapolation}
Next, we study mini-batch SGDE and its convergence. The updates of mini-batch SGDE are presented in Algorithm~\ref{alg:SGDEmini}. 
%
The convergence result of mini-batch SGDE is given in Theorem~\ref{thm:SGDE}.
\begin{thm}\label{thm:SGDE}
Under Assumption~\ref{ass:1} (i) (iii) and (iv), let $\eta \leq \frac{1}{12L}$ and $\x_1=\z_0= \x_0$, then SGDE ensures that
\begin{align}\label{upp:bound:SGDE}
\min_{t\in\{1,\ldots, T\}}\E[\|\nabla f(\x_t)\|^2]&\leq \frac{3L\eta G^2}{2T} +  \frac{8(f(\x_0) - f(\x_*))}{\eta T} + \frac{72G^2}{m}- \frac{1}{\eta^2T}\sum_{t=0}^{T-1}\E[\|\x_{t+1} -\x_{t}\|^2].
\end{align}
where $\x_* = \arg\min_{\x\in\R^d}f(\x)$. Under an additional assumption that $f(\x_0) - f(\x_*) \leq \Delta_0$ where $\Delta_0>0$ is a constant, in order to have $\min_{t\in\{1,\ldots, T\}}\E[\|\nabla f(\x_t)\|] \leq \epsilon$, the iteration complexity is $T=O(1/\epsilon^2)$ with mini-batch size $m=O(1/\epsilon^2)$, indicating that the gradient complexity is $O(1/\epsilon^4)$.
\end{thm}
{\bf Remark. }The gradient complexity $O(1/\epsilon^4)$ of mini-batch SGDE matches that of SGD method for stochastic non-convex optimization. However, comparing with the convergence upper bound of SGD, the above bound of GDE in (\ref{upp:bound:SGDE}) has an additional negative term $- \frac{1}{\eta^2T}\sum_{t=0}^{T-1}\E[\|\x_{t+1} -\x_{t}\|^2]$.
\begin{proof}
By the $L$-smooth of $f(\x)$ we have
\begin{align*}
f(\x_t)&\leq f(\z_{t-1}) + \nabla f(\z_{t-1})^{\top}(\x_t - \z_{t-1}) + \frac{L}{2}\|\x_t - \z_{t-1}\|^2\\
&\leq f(\z_{t-1}) + \nabla f(\x_t)^{\top}(\x_t - \z_{t-1}) + (\nabla f(\z_{t-1}) - \nabla f(\x_t))^{\top}(\x_t - \z_{t-1}) + \frac{L}{2}\|\x_t - \z_{t-1}\|^2\\
&\leq f(\z_{t-1}) + \nabla f(\x_t)^{\top}(\x_t - \z_{t-1}) + \frac{3 L}{2}\|\x_t - \z_{t-1}\|^2.
\end{align*}
Applying Lemma~\ref{lem:key1} with $\u=\z_{t-1}, \x=\x_t, \z=\z_{t-1}, \z_+=\z_t, \xi=\frac{1}{m}\sum_{i=1}^{m}\nabla f(\x_{t-1};\xi_{i,t-1}), \zeta=\frac{1}{m}\sum_{i=1}^{m}\nabla f(\x_t;\xi_{i,t}), \gamma=\eta$, we have
\begin{align*}
&\left(\frac{1}{m}\sum_{i=1}^{m}\nabla f(\x_t;\xi_{i,t})\right)^{\top}(\x_t - \z_{t-1})\\
\leq& -\frac{\|\z_t - \z_{t-1}\|^2}{2\eta} + \eta \bigg\|\frac{1}{m}\sum_{i=1}^{m}\nabla f(\x_t;\xi_{i,t}) - \frac{1}{m}\sum_{i=1}^{m}\nabla f(\x_{t-1};\xi_{i,t-1}) \bigg\|^2 \\
&- \frac{1}{2\eta}(\|\x_t - \z_{t-1}\|^2 + \|\x_t - \z_t\|^2).
\end{align*}
Taking expectation on both sides, we have
\begin{align*}
&\E\bigg[\left(\frac{1}{m}\sum_{i=1}^{m}\nabla f(\x_t;\xi_{i,t})\right)^{\top}(\x_t - \z_{t-1})\bigg]\\
\leq& \E\bigg[-\frac{\|\z_t - \z_{t-1}\|^2}{2\eta} + \eta \bigg\|\frac{1}{m}\sum_{i=1}^{m}\nabla f(\x_t;\xi_{i,t}) - \frac{1}{m}\sum_{i=1}^{m}\nabla f(\x_{t-1};\xi_{i,t-1}) \bigg\|^2\\
& - \frac{1}{2\eta}(\|\x_t - \z_{t-1}\|^2 + \|\x_t - \z_t\|^2)\bigg].
\end{align*}
Combining the above two inequalities together, we have
\begin{align*}
\E[f(\x_t)] \leq \E\bigg[ &f(\z_{t-1}) -\frac{\|\z_t - \z_{t-1}\|^2}{2\eta} + \eta \bigg\|\frac{1}{m}\sum_{i=1}^{m}\nabla f(\x_t;\xi_{i,t}) - \frac{1}{m}\sum_{i=1}^{m}\nabla f(\x_{t-1};\xi_{i,t-1}) \bigg\|^2  \\
&- \frac{1}{2\eta}(\|\x_t - \z_{t-1}\|^2 + \|\x_t - \z_t\|^2) + \frac{3L}{2}\|\x_t -\z_{t-1}\|^2\bigg].
\end{align*}
Moreover, by the smoothness of $f(\x)$
\begin{align*}
f(\z_t)&\leq f(\x_t) + \nabla f(\x_t)^{\top}(\z_t - \x_t) + \frac{L}{2}\|\x_t -\z_t\|^2\\
&\leq f(\x_t) +  \frac{\eta \|\nabla f(\x_t)\|^2}{4} + \frac{\|\x_t - \z_t\|^2}{\eta} + \frac{L}{2}\|\x_t - \z_t\|^2
\end{align*}
Combining the above two inequalities together and taking expectation on both sides, then
\begin{align*}
&\E[f(\z_t)- f(\z_{t-1})]\\
\leq &\E\bigg[ \frac{\eta\|\nabla f(\x_t)\|^2 }{4}  -\frac{\|\z_t - \z_{t-1}\|^2}{2\eta} +(\frac{1}{\eta} + \frac{L}{2})\|\x_t - \z_t\|^2 + \frac{3L}{2}\|\x_t -\z_{t-1}\|^2\bigg]\\
&+ \E\bigg[ \eta \bigg\|\frac{1}{m}\sum_{i=1}^{m}\nabla f(\x_t;\xi_{i,t}) - \frac{1}{m}\sum_{i=1}^{m}\nabla f(\x_{t-1};\xi_{i,t-1}) \bigg\|^2 - \frac{1}{2\eta}(\|\x_t - \z_{t-1}\|^2 + \|\x_t - \z_t\|^2) \bigg]\\
= &\E\left[\frac{\eta\|\nabla f(\x_t)\|^2 }{4}   -\frac{\eta\|\frac{1}{m}\sum_{i=1}^{m}\nabla f(\x_t;\xi_{i,t}) \|^2}{2}\right] + \frac{3\eta^2 L}{2}\E\bigg[\bigg\|\frac{1}{m}\sum_{i=1}^{m}\nabla f(\x_{t-1};\xi_{i,t-1}) \bigg\|^2\bigg] \\
& + \E\bigg[ \bigg(\eta+\frac{\eta^2 L}{2}+\eta\bigg)\bigg\|\frac{1}{m}\sum_{i=1}^{m}\nabla f(\x_t;\xi_{i,t}) - \frac{1}{m}\sum_{i=1}^{m}\nabla f(\x_{t-1};\xi_{i,t-1}) \bigg\|^2 \\
& - \frac{1}{2\eta}(\|\x_t - \z_{t-1}\|^2 + \|\x_t - \z_t\|^2)\bigg]
\end{align*}
where we use the facts that $\x_t - \z_{t-1} =-  \frac{\eta}{m}\sum_{i=1}^{m}\nabla f(\x_{t-1};\xi_{i,t-1}) $, $\z_t - \z_{t-1} = - \frac{\eta}{m}\sum_{i=1}^{m}\nabla f(\x_t;\xi_{i,t})$, and $\x_t -\z_t =  \frac{\eta}{m}\sum_{i=1}^{m}\nabla f(\x_t;\xi_{i,t}) - \frac{\eta}{m}\sum_{i=1}^{m}\nabla f(\x_{t-1};\xi_{i,t-1})$.
By $\E\bigg[\bigg \| \frac{1}{m}\sum_{i=1}^{m}\nabla f(\x;\xi_{i}) \bigg\|^2\bigg] = \| f(\x)\|^2 + \E\bigg[\bigg\|\frac{1}{m}\sum_{i=1}^{m}\nabla f(\x;\xi_{i})-f(\x)\bigg\|^2\bigg]$, and the assumption of $\E\bigg[\bigg\|\frac{1}{m}\sum_{i=1}^{m}\nabla f(\x;\xi_{i})-f(\x)\bigg\|^2\bigg] \leq \frac{G^2}{m}$, then
\begin{align*}
&\E[f(\z_t)- f(\z_{t-1})]\\
\leq &\E\left[ -\frac{\eta\|\nabla f(\x_t)\|^2 }{4}   -\frac{\eta\|\frac{1}{m}\sum_{i=1}^{m}\nabla f(\x_t;\xi_{i,t})-\nabla f(\x_t)\|^2}{2}\right] \\
&+ \frac{3\eta^2 L}{2}\E[\|\nabla f(\x_{t-1})\|^2]  + \frac{3\eta^2 L}{2}\E\bigg[\bigg\|\frac{1}{m}\sum_{i=1}^{m}\nabla f(\x_{t-1};\xi_{i,t-1})-\nabla f(\x_{t-1})\bigg\|^2\bigg] \\
&+ \E\bigg[ \bigg(2\eta+\frac{\eta^2 L}{2}\bigg)\bigg \|\frac{1}{m}\sum_{i=1}^{m}\nabla f(\x_t;\xi_{i,t}) -\nabla f(\x_t)- \frac{1}{m}\sum_{i=1}^{m}\nabla f(\x_{t-1};\xi_{i,t-1})+\nabla f(\x_{t-1})\bigg\|^2  \bigg]\\
& + \E\bigg[ \bigg(2\eta+\frac{\eta^2 L}{2}\bigg) \|\nabla f(\x_t) - \nabla f(\x_{t-1})\|^2  \bigg] - \frac{1}{2\eta} \E\bigg[ \|\x_t - \z_{t-1}\|^2 + \|\x_t - \z_t\|^2\bigg]\\
\leq &\E\left[ -\frac{\eta\|\nabla f(\x_t)\|^2 }{4}  \right] + \frac{3\eta^2 L}{2}\E[\|\nabla f(\x_{t-1})\|^2]  +\bigg(8\eta+ \frac{7\eta^2 L}{2}\bigg)\frac{G^2}{m} \\
& + \E\bigg[ \bigg(2\eta+\frac{\eta^2 L}{2}\bigg) \|\nabla f(\x_t) - \nabla f(\x_{t-1})\|^2  \bigg] - \frac{1}{2\eta} \E\bigg[ \|\x_t - \z_{t-1}\|^2 + \|\x_t - \z_t\|^2\bigg].
\end{align*}
Taking summation on both sides, we have 
\begin{align*}
& \E\bigg[\sum_{t=1}^T\frac{\eta}{4}\|\nabla f(\x_t)\|^2 - \frac{3L\eta^2}{2}\sum_{t=1}^T\|\nabla f(\x_{t-1})\|^2\bigg]\\
\leq& (2\eta + \frac{\eta^2 L}{2} )\sum_{t=0}^{T-1} \E\bigg[\|\nabla f(\x_{t+1}) - \nabla f(\x_{t})\|^2\bigg] +\bigg(8\eta+ \frac{7\eta^2 L}{2}\bigg)\frac{G^2T}{m}\\
& - \frac{1}{2\eta}\sum_{t=1}^T\E\bigg[ \|\x_t - \z_{t-1}\|^2 + \|\x_t - \z_t\|^2\bigg] +\sum_{t=1}^T \E[f(\z_{t-1}) - f(\z_t)].
\end{align*}
Note that 
\begin{align*}
&\sum_{t=1}^T(\|\x_t - \z_{t-1}\|^2  + \|\x_t - \z_t\|^2)\\
=&\sum_{t=1}^{T-1}\|\x_{t+1} - \z_{t}\|^2  + \|\x_{1} - \z_{0}\|^2 + \sum_{t=1}^T\|\x_t - \z_t\|^2 \\
\geq& \frac{1}{2}\sum_{t=1}^{T-1}\|\x_t - \x_{t+1}\|^2+ \|\x_{1} - \z_{0}\|^2\\
=  &\frac{1}{2}\sum_{t=0}^{T-1}\|\x_t - \x_{t+1}\|^2 + \frac{1}{2}\|\x_1 - \x_{0}\|^2.~~(\text{Since }\z_0 = \x_0)
\end{align*}
Combining the above inequalities together, we have
\begin{align*}
&\sum_{t=1}^T\bigg(\frac{\eta}{4} - \frac{3L\eta^2}{2}\bigg)\E[\|\nabla f(\x_t)\|^2]\\
\leq& \frac{3L\eta^2}{2}\|\nabla f(\x_0)\|^2 -  \frac{1}{4\eta}\|\x_1 - \x_{0}\|^2+\sum_{t=1}^T \E[f(\z_{t-1}) - f(\z_t)]\\
&+(2\eta L^2 + \frac{\eta^2 L^3}{2} - \frac{1}{4\eta})\sum_{t=0}^{T-1}\E[\|\x_{t+1} -\x_{t}\|^2]+\bigg(8\eta+ \frac{7\eta^2 L}{2}\bigg)\frac{G^2T}{m}.
\end{align*}
Since $\eta \leq \frac{1}{12L}$, then $8\eta+ \frac{7\eta^2 L}{2}\leq 9\eta$, $2\eta L^2  + \frac{\eta^2 L^3}{2}- \frac{1}{4\eta} \leq - \frac{1}{8\eta} $ and $\frac{\eta}{4} - \frac{3L\eta^2}{2} \geq \frac{\eta}{8}$.
Note that we can define $\nabla f(\x_0;\xi_0)=0$ (i.e, $\x_1 = \z_0$), which will not affect our analysis above. Then $\x_1=\z_0= \x_0$,  and 
\begin{align*}
\frac{\eta}{8}\sum_{t=1}^T\E\|\nabla f(\x_t)\|^2]&\leq \frac{3L\eta^2G^2}{2} + f(\z_0) - f(\x_*) + \frac{9G^2\eta T}{m}  - \frac{1}{8\eta}\sum_{t=0}^{T-1}\E[\|\x_{t+1} -\x_{t}\|^2] ,
\end{align*}
which implies
\begin{align*}
\min_{t\in\{1,\ldots, T\}}\E[\|\nabla f(\x_t)\|^2]&\leq \frac{3L\eta G^2}{2T} +  \frac{8(f(\z_0) - f(\x_*))}{\eta T} + \frac{72G^2}{m}- \frac{1}{\eta^2T}\sum_{t=0}^{T-1}\E[\|\x_{t+1} -\x_{t}\|^2].
\end{align*}
\end{proof}

\subsection{Stagewise SGDE}
In the previous subsection, mini-batch SGDE requires the mini-batch size in the order of $O({1}/{\epsilon^2})$, which might be not practical when the target accuracy $\epsilon$ is sufficiently small. 
In this subsection, we propose a new variant of SGDE without requiring a large mini-batch size, and we present the details in Algorithm~\ref{alg:stagewise:SGDE:2} with a subroutine SGDE in Algorithm~\ref{alg:SGDE2}, which is referred to stagewise SGDE. For $s$-th stage, stagewise SGDE solves the following subproblem approximately
\begin{align*}
f_s(\x) = f(\x) + \frac{1}{2\gamma}\|\x-\x^{s-1}\|^2,
\end{align*}
where $\x^{s-1}$ is the solution of the last stage, and $\gamma = \frac{1}{4L}$ is a constant. It is easy to show that $f_s(\x)$ is convex under the Assumption~\ref{ass:1} (i), meaning that one may employ SGDE algorithm with convergence guarantee for convex problems. By using the convexity of $f_s$, the subroutine SGDE usually returns an average solution. Besides, stagewise SGDE uses a decreasing sequence of step size $\eta_s$ and an increasing sequence of iteration number $T_s$. Different from GDE and mini-batch SGDE, the final  solution of stagewise SGDE is selected from the sequence of stagewise averaged solutions $\{\x_s\}$ based on non-uniform sampling probabilities increasing as the stage number $s$. It is notable that this type of stagewise algorithm has been investigated in existing studies~(see~\citep{chen2018universal} and references therein). However, to the best of our knowledge, the proposed algorithm is the first work that runs SGDE method in a stagewise manner with the theoretical guarantee for non-convex optimization. 
We present the convergence result of stagewise SGDE in Theorem~\ref{thm:stagewise:sgde:2}.
\begin{algorithm}[t]
\caption{SGDE($\x_0, f, \eta, T$)} \label{alg:SGDE2}
\begin{algorithmic}
\STATE Initilization: $\z_0 = \x_0$ and $\g_0 = \nabla f(\x_{0}; \xi_{0}) $
\FOR{$t=1,\ldots, T$}
    	\STATE $\x_t = \z_{t-1} - \eta \g_{t-1}$
	\STATE $\g_t = \nabla f(\x_{t}; \xi_{t}) $
	\STATE $\z_t =  \z_{t-1} - \eta \g_t$
\ENDFOR
\RETURN $\widehat \x_T = \frac{1}{T}\sum_{t=1}^{T} \x_t$
\end{algorithmic}
\end{algorithm}
\begin{algorithm}[t]
\caption{stagewise SGDE} \label{alg:stagewise:SGDE:2}
\begin{algorithmic}
\STATE Initilization: $\x^0 = \x_0$
\FOR{$s=1,\ldots, S$}
    	\STATE $f_s(\x) = f(\x) + \frac{1}{2\gamma}\|\x-\x^{s-1}\|^2$
	\STATE $\x^s  =  \text{SGDE}(\x^{s-1}, f_s, \eta_s, T_s)$
\ENDFOR
\STATE \textbf{Return:} $\x_{\tau}$, $\tau$ is randomly chosen from $\{1, \ldots, S\}$ according to probabilities $p_\tau = \frac{w_{\tau}}{\sum_{s=1}^{S}w_{s}}, \tau=1, \ldots, S$.
\end{algorithmic}
\end{algorithm}
\begin{thm}\label{thm:stagewise:sgde:2}
Under Assumption~\ref{ass:1}, by running Algorithm~\ref{alg:stagewise:SGDE:2} with $\gamma = \frac{1}{4L}$, $w_s = s^\alpha$ ($\alpha>1$),  $\eta_s = \frac{c \gamma}{3s} \leq \frac{1}{2L} = \frac{\gamma}{3}$, and $T_s = \frac{36s}{c}$, then
\begin{align}\label{upp:bound:StagewiseSGDE}
 \E[\|\nabla f(\x_{\tau})\|^2]\leq 
 \left\{\begin{array}{cc}\frac{20\Delta(\alpha+1)}{\gamma(S+1)}     +  \frac{480 G^2 c (\alpha+1)}{S+1} - \frac{ 60 \sum_{s=1}^{S+1}w_s D_{T_s}}{\gamma\sum_{s=1}^{S+1}w_s},& \alpha\geq 1,\\
 \frac{20\Delta(\alpha+1)}{\gamma(S+1)}    +   \frac{480 G^2 c (\alpha+1)}{\alpha  (S+1)} - \frac{ 60  \sum_{s=1}^{S+1} w_s D_{T_s}}{\gamma \sum_{s=1}^{S+1}w_s},& 0< \alpha < 1,
 \end{array}\right.
\end{align}
where $D_{T_s} = \frac{1}{16 T_s \eta_s}\sum_{t=1}^{T_s} \|\x_t - \x_{t-1}\|^2$.
Therefore, in order to have $\E[\|\nabla f(\x_{\tau})\|^2]\leq \epsilon^2$, we can set $S=O(1/\epsilon^2)$. The total number of iterations is $O\left(\frac{1}{\epsilon^4}\right)$.
\end{thm}
{\bf Remark. }Although the iteraction complexity of stagewise SGDE mathches that of stagewise SGD in~\citep{chen2018universal}, the above bound of stagewise SGDE in (\ref{upp:bound:StagewiseSGDE}) has an additional negative term $- \frac{ 60  \sum_{s=1}^{S+1} w_s D_{T_s}}{\gamma \sum_{s=1}^{S+1}w_s}$, comparing with the convergence upper bound of stagewise SGD. This negative term could help improve convergence in practice. 
\begin{proof}
For the $s$-th stage, the following problem is solved
\begin{align*}
\min_{\x} f_s(\x)=f(\x) + \frac{1}{2\gamma}\|\x - \x_{s-1}\|^2
\end{align*}
where $\x_{s-1}$ is the solution from last stage. Let define $\widehat\z_s = \arg\min_{\x} f_s(\x)$.
By applying Lemma~\ref{lem:key1} with $\u=\widehat\z_s, \x=\x_t, \z=\z_{t-1}, \z_+=\z_t, \xi=\nabla f_s(\x_{t-1}; \xi_{t-1}), \zeta=\nabla f_s(\x_t; \xi_t), \gamma=\eta_s$, we have
\begin{align*}
\nabla f_s(\x_t; \xi_t)^{\top}(\x_t - \widehat\z_s)\leq& \frac{\|\widehat\z_s - \z_{t-1}\|^2-\|\widehat\z_s - \z_{t}\|^2}{2\eta_s} + \eta_s \|\nabla f_s(\x_t; \xi_t) - \nabla f_s(\x_{t-1}; \xi_{t-1})\|^2 \\
&- \frac{1}{2\eta_s}(\|\x_t - \z_{t-1}\|^2 + \|\x_t - \z_t\|^2).
\end{align*}
Taking average over $t=1,\dots,T_s$ for above inequality and by the convexity of $f(\x)$ we have
\begin{align*}
&\frac{1}{T_s}\sum_{t=1}^{T_s}\nabla f_s(\x_t; \xi_t)^{\top}(\x_t - \widehat\z_s)
\leq \frac{\|\widehat\z_s - \z_{0}\|^2}{2\eta_sT_s} - \frac{1}{2\eta_sT_s}\sum_{t=1}^{T_s}( \|\x_t - \z_{t-1}\|^2 + \|\x_t - \z_t\|^2)\\
&+ \frac{\eta_s}{T_s}\sum_{t=1}^{T_s} \|\nabla f(\x_t; \xi_t) - \nabla f(\x_{t-1}; \xi_{t-1}) + \frac{1}{\gamma}(\x_t - \x_{t-1})\|^2\\
\leq& \frac{\|\widehat\z_s - \z_{0}\|^2}{2\eta_sT_s} - \frac{1}{2\eta_sT_s}\sum_{t=1}^{T_s}( \|\x_t - \z_{t-1}\|^2 + \|\x_t - \z_t\|^2)\\
&+ \frac{2\eta_s}{T_s}\sum_{t=1}^{T_s} \|\nabla f(\x_t; \xi_t) - \nabla f(\x_{t-1}; \xi_{t-1})\|^2 + \frac{2\eta_s}{\gamma^2T_s}\sum_{t=1}^{T_s} \|\x_t - \x_{t-1}\|^2\\
\leq &  \frac{\|\widehat\z_s - \z_0\|^2}{2\eta_sT_s} - \frac{1}{2\eta_sT_s}\sum_{t=1}^{T_s}(\|\x_t - \z_{t-1}\|^2 + \|\x_t - \z_t\|^2) +  \frac{6\eta_s}{T_s}\sum_{t=1}^{T_s} \|\nabla f(\x_t) - \nabla f(\x_{t-1})\|^2\\
&+  \frac{6\eta_s}{T_s}\sum_{t=1}^{T_s}( \|\nabla f(\x_t; \xi_t) - \nabla f(\x_t)\|^2 +  \|\nabla f(\x_{t-1}; \xi_{t-1}) - \nabla f(\x_{t-1})\|^2)+ \frac{2\eta_s}{\gamma^2T_s}\sum_{t=1}^{T_s} \|\x_t - \x_{t-1}\|^2,
\end{align*}
where the last inequality is due to $\| \x_1 + \x_2 + \x_3\|^2_2 \leq 3(\|\x_1\|_2^2+\|\x_2\|_2^2+\|\x_3\|_2^2)$.
By the smoothness of $f(\x)$ we have $\frac{\eta_s}{T_s}\sum_{t=1}^{T_s} \|\nabla f(\x_t) - \nabla f(\x_{t-1})\|^2 \leq \frac{\eta_sL^2}{T_s}\sum_{t=1}^{T_s} \|\x_t - \x_{t-1}\|^2$. Similarly, note that $\sum_{t=1}^{T_s} (\|\x_t - \z_{t-1}\|^2  + \|\x_t - \z_t\|^2)\geq\frac{1}{2}\sum_{t=1}^{T_s}\|\x_t - \x_{t-1}\|^2 +  \frac{1}{2}\|\x_1 - \x_{0}\|^2 \geq \frac{1}{2}\sum_{t=1}^{T_s}\|\x_t - \x_{t-1}\|^2$. 
Let $\Delta_t := \nabla f(\x_t; \xi_t) - \nabla f(\x_t)$, and by setting of $\gamma = 1/(4L)$, then the above iequality becomes
\begin{align}\label{res1:stage:ext}
\nonumber & \frac{1}{T_s}\sum_{t=1}^{T_s}\nabla f_s(\x_t; \xi_t)^{\top}(\x_t - \widehat\z_s)\\
\nonumber  \leq &  \frac{\|\widehat\z_s - \z_0\|^2}{2\eta_sT_s} +  \frac{38\eta_s L^2 - \frac{1}{4\eta_s}}{T_s}\sum_{t=1}^{T_s} \|\x_t - \x_{t-1}\|^2 + \frac{6\eta_s}{T_s}\sum_{t=1}^{T_s}( \|\Delta_t\|^2 +  \|\Delta_{t-1}\|^2)\\
 \leq &  \frac{\|\widehat\z_s - \z_0\|^2}{2\eta_sT_s} + \frac{6\eta_s}{T_s}\sum_{t=1}^{T_s}( \|\Delta_t\|^2 +  \|\Delta_{t-1}\|^2) - \underbrace{\frac{1}{16 T_s \eta_s}\sum_{t=1}^{T_s} \|\x_t - \x_{t-1}\|^2}\limits_{D_{T_s}} ,
\end{align}
where the last inequality is due to $\eta_s \leq \frac{1}{16L}$ so that $3\eta_s L^2 - \frac{1}{4\eta_s} \leq - \frac{1}{16\eta_s} $. 
Since $\E[\nabla f_s(\x_t; \xi_t)^{\top}(\x_t - \widehat\z_s) | \x_t, \Delta_{t-1}, \dots, \Delta_0] = \nabla f_s(\x_t)^{\top}(\x_t - \widehat\z_s)$, and $\E[\|\Delta_t\|^2 | \x_t, \Delta_{t-1}, \dots, \Delta_0] \leq G^2$, then by the convexity of $f_s(\x)$ we have
\begin{align*}
&\E\bigg[f_s(\x_s) - f_s(\widehat\z_s)\bigg]\leq \E\bigg[\frac{1}{T_s}\sum_{t=1}^{T_s}f_s(\x_t) - f_s(\widehat\z_s)\bigg]\\
 \leq&  \E\bigg[\frac{1}{T_s}\sum_{t=1}^{T_s}\nabla f_s(\x_t)^{\top}(\x_t - \widehat\z_s)\bigg] =  \E\bigg[\frac{1}{T_s}\sum_{t=1}^{T_s}\E[\nabla f_s(\x_t; \xi_t)^{\top}(\x_t - \widehat\z_s) | \x_t, \Delta_{t-1}, \dots, \Delta_0]\bigg]\\
\leq &\frac{\E[\|\widehat\z_s - \z_0\|^2]}{2\eta_sT_s} + 12\eta_sG^2 - D_{T_s}= \frac{\E[\|\widehat\z_s - \x_{s-1}\|^2]}{2\eta_sT_s} + 12\eta_sG^2- D_{T_s},
\end{align*}
Since $f_s(\x_s) =  f(\x_s) + \frac{1}{2\gamma}\|\x_s - \x_{s-1}\|^2 $ and $f_s(\widehat\z_s) \leq f(\x_{s-1})$, then
\begin{align*}
\E\bigg[f(\x_s) + \frac{1}{2\gamma}\|\x_s - \x_{s-1}\|^2 - f(\x_{s-1})\bigg] \leq& \frac{\E[\|\widehat\z_s - \x_{s-1}\|^2]}{2\eta_sT_s} + 12\eta_sG^2- D_{T_s}.
\end{align*}
By Young's inequality $\|\x_s - \x_{s-1}\|^2 \geq \frac{1}{2} \|\widehat\z_s - \x_{s-1}\|^2 - \|\x_s -\widehat\z_s\|^2$, then
\begin{align*}
&\E\bigg[ \bigg( \frac{1}{4\gamma} - \frac{1}{2\eta_s T_s}\bigg)\|\widehat\z_s - \x_{s-1}\|^2  \bigg] \\
\leq& \frac{1}{2\gamma}\E[\|\x_s -\widehat\z_s\|^2]+ 12\eta_sG^2 + \E\bigg[ f(\x_{s-1}) - f(\x_s)\bigg]- D_{T_s}\\
\leq& \frac{1}{\gamma(\gamma^{-1} - \mu)}\E[f_s(\x_s) -f_s(\widehat\z_s)]+ 12\eta_sG^2 + \E\bigg[ f(\x_{s-1}) - f(\x_s)\bigg]- D_{T_s}\\
\leq& \frac{1}{\gamma(\gamma^{-1} - \mu)}\bigg(\frac{\E[\|\widehat\z_s - \x_{s-1}\|^2]}{2\eta_sT_s} + 12\eta_sG^2 - D_{T_s} \bigg) + 12\eta_sG^2 + \E\bigg[ f(\x_{s-1}) - f(\x_s)\bigg]- D_{T_s},
\end{align*}
where the second inequality uses the $(\gamma^{-1}-\mu)$-strong convex of $f_s(\x)$ and the last inequality uses (\ref{res1:stage:ext}).
By setting $\gamma^{-1} = 2\mu$, then the above inequality will be
\begin{align*}
\E\bigg[ \bigg( \frac{1}{4\gamma} - \frac{3}{2\eta_s T_s}\bigg)\|\widehat\z_s - \x_{s-1}\|^2  \bigg] \leq 36\eta_sG^2 + \E\bigg[ f(\x_{s-1}) - f(\x_s)\bigg]- 3D_{T_s},
\end{align*}
As long as $\eta_sT_s\geq 12\gamma$, we have
\begin{align*}
\frac{1}{8\gamma}\E\bigg[ \|\widehat\z_s - \x_{s-1}\|^2  \bigg] \leq 36\eta_sG^2 + \E\bigg[ f(\x_{s-1}) - f(\x_s)\bigg]-3D_{T_s}.
\end{align*}
By the property of  Moreau envelope funtion, we know $\nabla f_{\gamma}(\x_{s-1}) =\frac{1}{\gamma}  (\widehat\z_s - \x_{s-1})$, then
\begin{align*}
\frac{\gamma}{8}\E\bigg[ \|\nabla f_{\gamma}(\x_{s-1})\|^2  \bigg] \leq 36\eta_sG^2 + \E\bigg[ f(\x_{s-1}) - f(\x_s)\bigg]- 3D_{T_s}.
\end{align*}
which implies
\begin{align*}
 \E\left[w_s \|\nabla f_{\gamma}(\x_{s-1}) \|^2  \right]   \leq  16\mu \E\left[ w_s(f(\x_{s-1}) - f(\x_s) )\right] + 576\mu w_s \eta_sG^2 - 48\mu w_s D_{T_s}.
\end{align*} 
By summing over $s=1,\dots, S+1$ we get 
\begin{align*}
 &\sum_{s=1}^{S+1}\E\left[w_s \|\nabla f_{\gamma}(\x_{s-1}) \|^2  \right]  \\
  \leq & 16\mu \E\left[  \sum_{s=1}^{S+1}w_s(f(\x_{s-1}) - f(\x_s) )\right] +  \sum_{s=1}^{S+1}576\mu w_s \eta_sG^2 - \sum_{s=1}^{S+1}48\mu w_s D_{T_s}.
\end{align*} 
Then taking the expectation over $\tau$, it becomes
\begin{align*}
 \E\left[\|\nabla f_{\gamma}(\x_{\tau}) \|^2  \right]   \leq  16\mu \E\left[  \frac{\sum_{s=1}^{S+1}w_s(f(\x_{s-1}) - f(\x_s) )}{\sum_{s=1}^{S+1}w_s}\right] +  \frac{\sum_{s=1}^{S+1}576\mu w_s\eta_s G^2}{\sum_{s=1}^{S+1}w_s} - \frac{ \sum_{s=1}^{S+1}48\mu w_s D_{T_s}}{\sum_{s=1}^{S+1}w_s}.
\end{align*} 
Let consider the term $\sum_{s=1}^{S+1}w_s(f(\x^{s-1}) - f(\x^s) )$:
\begin{align*}
&\sum_{s=1}^{S+1}w_s(f(\x^{s-1}) - f(\x^s) ) \\
=& \sum_{s=1}^{S+1}[w_{s-1}f(\x^{s-1}) - w_sf(\x^s)] + \sum_{s=1}^{S+1}(w_s-w_{s-1})f(\x^{s-1})\\
=& w_{0}f(\x^{0}) - w_{S+1}f(\x^{S+1}) + \sum_{s=1}^{S+1}(w_s-w_{s-1})f(\x^{s-1})\\
=& w_{0}(f(\x^{0}) -f(\x_*) )- w_{S+1}(f(\x^{S+1})-f(\x_*) ) + \sum_{s=1}^{S+1}(w_s-w_{s-1})(f(\x^{s-1})-f(\x_*) )\\
\leq & w_0 \Delta + 0 + \sum_{s=1}^{S+1}(w_s-w_{s-1})\Delta = w_{S+1}\Delta
\end{align*}
Then,
\begin{align*}
 \E\left[\|\nabla f_{\gamma}(\x^{\tau}) \|^2  \right]   \leq   \frac{16\mu  w_{S+1}\Delta}{\sum_{s=1}^{S+1}w_s} +  \frac{\sum_{s=1}^{S+1}576\mu w_s\eta_s G^2}{\sum_{s=1}^{S+1}w_s} - \frac{ \sum_{s=1}^{S+1}48\mu w_s D_{T_s}}{\sum_{s=1}^{S+1}w_s}.
\end{align*} 
We know $w_s = s^\alpha$ $(\alpha >1)$, the standard calculus tells
\begin{align*}
&\sum_{s=1}^Ss^{\alpha}\geq \int_0^{S}x^\alpha d x= \frac{S^{\alpha+1}}{\alpha+1},~\forall\alpha>0,\\
&\sum_{s=1}^Ss^{\alpha-1}\leq S^\alpha,~\forall\alpha\geq 1,\\
&\sum_{s=1}^Ss^{\alpha-1}\leq \int_{0}^S x^{\alpha-1}d x = \frac{ S^\alpha}{\alpha},~\forall 0<\alpha<1.
\end{align*}
Since $\eta_s = \frac{c}{Ls} < \frac{1}{2L}$, $(L=3\mu = \frac{3}{2\gamma})$ then
\begin{align*}
 \E[\|\nabla f_{\gamma}(\x^{\tau})\|^2]\leq 
 \left\{\begin{array}{cc}\frac{8\Delta(\alpha+1)}{\gamma(S+1)}     +  \frac{192 G^2 c (\alpha+1)}{S+1}- \frac{ \sum_{s=1}^{S+1}48\mu w_s D_{T_s}}{\sum_{s=1}^{S+1}w_s}& \alpha\geq 1,\\
 \frac{8\Delta(\alpha+1)}{\gamma(S+1)}    +   \frac{192 G^2 c (\alpha+1)}{\alpha (S+1)}- \frac{ \sum_{s=1}^{S+1}48\mu w_s D_{T_s}}{\sum_{s=1}^{S+1}w_s}& 0< \alpha < 1.
 \end{array}\right.
\end{align*}
By the result of Moreau envelop function in (\ref{moreau:ineqs}), we know for any $\x$
\begin{align*}
\|\nabla f (\x)\| \leq & \|\nabla f (\x) - \nabla f (\widehat \x) \| + \|\nabla f(\widehat \x)\| \\
\leq & L \| \x - \widehat \x\| + \|\nabla f_\gamma(\x)\| \\
= & (1+ L \gamma) \|\nabla f_\gamma(\x)\| = \frac{5}{2}\|\nabla f_\gamma(\x)\|.
\end{align*}
Therefore, in order to have $\E[\|\nabla f(\x^{\tau})\|^2]\leq \epsilon^2$, i.e., $\E[\|\nabla f_{\gamma}(\x^{\tau})\|^2]\leq \frac{4}{25} \epsilon^2$, we can set $S=O(1/\epsilon^2)$. The total number of iterations is
\begin{align*}
\sum_{s=1}^{S} T_s = \sum_{s=1}^{S}\frac{36 s}{c} = O\left(\frac{1}{\epsilon^4}\right).
\end{align*}
\end{proof}

\section{Conclusions}
In this paper, we have presented a GDE algorithm for solving smooth non-convex optimization problems and two stochastic variants of mini-batch GDE namely SGDE and stagewise SGDE for solving smooth non-convex stochastic optimization problems. We have established their convergence results in terms of finding an approximate first-order stationary point. In particular, we provided convergence upper bounds of the proposed algorithms as the theoretical evidence on the advantage of extrapolation steps.

%


\bibliography{all-extragrad}
\end{document}